\documentclass [12pt,a4paper]{article}
\usepackage{amssymb, theorem}
\usepackage{graphicx}
\usepackage{graphics}
\newtheorem{Thm}{Theorem}
\newtheorem{prop}{Proposition}
\newtheorem{lemma}{Lemma}
\theorembodyfont{\rm}

\def\proof{{\it Proof: }}

\def\qed{\nobreak\hfill $\square$}
\def\osum{\oplus}
\def\bal{\langle\!\langle}
\def\jobb{\rangle\!\rangle}
\def\<{\langle}
\def\>{\rangle}
\def\aa{\alpha}

\def\ffi{\varphi}
\def\Cov{\mathrm{Cov}}
\def\Ima{\mathrm{Im}\,}
\def\Det{\mathrm{Det}}
\def\QCov{\mathrm{qCov}}
\def\pard{\partial}

\def\iH{{\cal H}}

\def\iP{{\cal P}}

\def\iM{{\cal M}}
\def\iN{{\cal N}}

\def\im{\mathrm{i}}

\def\bbbr{{\mathbb R}}

\def\Diag{\mbox{Diag}\,}

\def\Tr{\mathrm{Tr}\,}
\def\pont{{\, \cdot \,}}
\begin{document}
%\rightline{\today }
 \ \vskip 1cm 
\centerline{\LARGE {\bf From quasi-entropy to skew information}}
\bigskip
\bigskip\bigskip
\bigskip
\centerline{\large D\'enes Petz\footnote{E-mail: petz@math.bme.hu.
Partially supported by the Hungarian Research Grant OTKA  T068258.}$^{,3,4}$ and
V.E. S\'andor Szab\'o\footnote{E-mail: sszabo@math.bme.hu.
Partially supported by the Hungarian Research Grant OTKA TS-49835.}$^{,4}$}
\bigskip
\begin{center}
$^3$ Alfr\'ed R\'enyi Institute of Mathematics, \\H-1364 Budapest,
POB 127, Hungary
\end{center}
\medskip
\centerline{$^4$ Department for Mathematical Analysis, BUTE,}
\centerline{H-1521 Budapest, POB 91, Hungary}
\bigskip
\bigskip\bigskip
\begin{abstract}
This paper gives an overview about particular quasi-entropies, generalized 
quantum covariances, quantum Fisher informations, skew-informations and their
relations. The point is the dependence on operator monotone functions. It 
is proven that a skew-information is the Hessian of a quasi-entropy.
The skew-information and some inequalities are extended to a von Neumann 
algebra setting. 
\smallskip


\end{abstract}

\bigskip\bigskip

\section{Introductory preliminaries}
Let $\iM$ denote the algebra of $n \times n$ matrices with complex entries.
For positive definite matrices $D_1, D_2\in \iM$, for $A \in \iM$ and 
a function $F:\bbbr^+ \to \bbbr$, the {\bf quasi-entropy} is defined as
\begin{eqnarray}\label{E:quasi}
S^A_F (D_1,D_2)&:=& \< A D_1^{1/2}, F(\Delta(D_2/ D_1))(AD_1^{1/2})\> \cr
&=&\Tr  D_1^{1/2} A^*F(\Delta(D_2/ D_1))(AD_1^{1/2}),
\end{eqnarray}
where $\Delta(D_2/ D_1):\iM \to \iM$ is a linear mapping acting 
on matrices:
$$
\Delta(D_2/ D_1)A=D_2 AD_1^{-1}.
$$
This concept was introduced in \cite{PD26, PD32}, see also Chapter 7 
in \cite{OP}
and it is the quantum generalization of the $F$-entropy of Csisz\'ar used 
in classical information theory (and statistics) \cite{Csi, LV}.
 
The concept of quasi-entropy includes some important special cases. If 
$D_1$ and $D_2$ are different and $A=I$, then we have a kind of relative 
entropy. For $F(x)=- \log x$ we have Umegaki's relative entropy 
$S(D_1\|D_2)=\Tr D_1 (\log D_1 - \log D_2)$. More generally,
$$
F(x)={1 \over \aa(1-\aa)}\big(1-x^{\aa}\big),
$$ 
is operator monotone decreasing for $\alpha \in (-1,1)$. (For $\aa=0$, the 
limit is taken and it is $- \log x$.) Then the R\'enyi entropies are 
produced 
$$
S_\alpha(D_1\|D_2):={1 \over \aa(1-\aa)}
\Tr (I-D_2^{\aa}D_1^{-\aa})D_1.
$$

If $D_1=D_2=D$ and $A, B \in \iM$ are arbitrary, then one can approach  to 
the {\bf generalized covariance} \cite{PD22}. An operator monotone function 
$f:\bbbr^+ \to \bbbr^+$ will be called {\bf standard} if $xf(x^{-1})=f(x)$ 
and $f(1)=1$. A standard function $f$ admits a canonical representation
\begin{equation}\label{E:canonicalf}
    f(t)=e^\beta\frac{1+t}{\sqrt{2}}\exp\int_0^1\frac{\lambda^2-1}{\lambda^2+1}\cdot
    \frac{1+t^2}{(\lambda+t)(1+\lambda t)}h(\lambda)\,d\lambda,
\end{equation}
where $ h:[0,1]\to[0,1] $ is a measurable function and $ \beta $ is a real 
constant \cite{H1}.

If $f$ is a standard function, then 
\begin{equation}\label{E:qC}
\QCov^f_{D}(A,B):=\< A D^{1/2}, f(\Delta(D/ D))(BD^{1/2})\>
-(\Tr DA^*)(\Tr DB).
\end{equation}
is a generalized covariance. The usual {\bf symmetrized covariance} 
corresponds to the function $f(t)=(t+1)/2$:
$$
\Cov_D(A,B):=
\frac{1}{2}\Tr (D(A^*B+BA^*))- (\Tr DA^*)(\Tr DB).
$$

The {\bf quantum Fisher information} is similarly defined to (\ref{E:quasi}),
but $F(x)=1/f(x)$ for a standard function $f:\bbbr^+ \to \bbbr^+$:
\begin{equation}\label{E:qF}
\gamma_D^f (A, B)=\< A D^{-1/2}, \frac{1}{f}(\Delta(D/ D))(BD^{-1/2})\>
\end{equation}
Quantum Fisher information was characterized by the monotonicity under 
coarse-graining \cite{PD2}. This kind of non-affine parametrization 
was used in \cite{PD2, PD22}, since the relation to operator means was 
emphasized.  Sometimes the affine parametrization is more convenient
and {\bf Hansen's canonical representation} of the inverse
of a standard operator monotone function can be used \cite{H}.

\begin{prop}\label{T:hansen}
If $f:\bbbr^+ \to \bbbr^+$ be a standard operator monotone function, then
\begin{equation}\label{E:canonical}
\frac{1}{f(t)}=
\int_0^1 \frac{1+\lambda}{2}
\left(\frac{1}{t+\lambda}+\frac{1}{1+t\lambda}\right)d\mu(\lambda),
\end{equation}
where $ \mu $ is a probability measure on $[0,1]$.
\end{prop}

The theorem implies that the set $\{1/f: f \mbox{\ is\ standard\ operator
\ monotone}\}$ is convex and gives the extremal points
\begin{equation}\label{E:efl}
g_\lambda(x):=\frac{1+\lambda}{2}
\left(\frac{1}{t+\lambda}+\frac{1}{1+t\lambda}\right) 
\qquad (0 \le \lambda \le 1).
\end{equation}
One can compute directly that
$$
\frac{\partial}{\partial\lambda}g_\lambda(x)=
-\frac{(1-\lambda^{2})(x+1)(x-1)^{2}}{2(x+\lambda)^{2}(1+x\lambda)^{2}}.
$$
Hence $g_\lambda$ is decreasing in the parameter $\lambda$. For $\lambda=0$
we have the largest function $g_0(t)=(t+1)/(2t)$ and for $\lambda=1$ the 
smallest is $g_1(t)=2/(t+1)$. Note that this was also obtained in the 
setting of positive operator means \cite{Ando}, harmonic and arithmetic means.

Covariance and Fisher information are bilinear (or sesqui-linear) forms, in
the applications they are mostly reduced to self-adjoint matrices.

The space $\iM$ has an orthogonal decomposition 
$$
\{B\in \iM : [D,B]=0\}\oplus \{\im[D,A]:A \in \iM\}.
$$
We denote the two subspaces by $\iM_D$ and $\iM_D^c$, respectively. If
$A_2 \in \iM_D$, then 
$$
F(\Delta(D/ D))(A_2D^{\pm 1/2})=A_2D^{\pm 1/2}
$$
implies
$$
\QCov^f_{D}(A_1,A_2)=\Tr DA_1^*A_2 -(\Tr DA_1^*)(\Tr DA_2),
\qquad
\gamma_D^f (A_1,A_2)=\Tr D^{-1}A_1^*A_2
$$
independently of the function $f$. Moreover, if $A_1 \in \iM_D^c$, then
$$
\gamma_D^f (A_1,A_2)=\QCov^f_{D}(A_1,A_2) =0 \,.
$$
Therefore, the effect of the function $f$ and the really quantum situation
are provided by the components from $\iM_D^c$.

\section{Quasi-entropy}

The quasi-entropies are monotone and jointly convex \cite{OP, PD32}:

Let $\alpha :\iM_0 \to \iM$ be a mapping between two matrix algebras. The dual
$\alpha^*: \iM \to \iM_0$ with respect to the Hilbert-Schmidt inner product
is positive if and only if $\alpha$ is positive. Moreover, $\alpha$ is
unital if and only if $\alpha^*$ is trace preserving. $\alpha: \iM_0 \to \iM$
is called a {\bf Schwarz mapping} if 
\begin{equation}\label{E:Sch}
\alpha(B^*B)\ge \alpha(B^*)\alpha(B)
\end{equation}
for every $B \in \iM_0$. 

\begin{prop}\label{P:quasimon}
Assume that $F:\bbbr^+ \to \bbbr$ is an operator monotone function
with $F(0)\ge 0$ and $\alpha:\iM_0 \to \iM$ is a unital Schwarz mapping. 
Then
\begin{equation}\label{E:quasimon}
S^A_F (\alpha^*(D_1),\alpha^*(D_2)) \ge
S^{\alpha(A)}_F (D_1,D_2)
\end{equation}
holds for $A \in \iM_0$ and for invertible density matrices $D_1$ and $D_2$
from the matrix algebra $\iM$.
\end{prop}

If we apply the monotonicity (\ref{E:quasimon}) to the embedding
$\alpha(X)=X \osum X$ of $\iM$ into $\iM \osum \iM$ and to the
densities $D_1=\lambda E_1\osum(1-\lambda)F_1$, $D_2=\lambda E_2\osum (1-\lambda)F_2$, 
then we obtain the joint concavity of the quasi entropy.

\begin{prop}\label{P:quasiconc}
Under the conditions of Theorem \ref{P:quasimon}, the joint concavity
\begin{equation}\label{E:quasiconc}
\lambda S^A_F (E_1,E_2)+(1-\lambda) S^A_F (F_1,F_2) 
\le S^{A}_F (\lambda E_1 +(1 -\lambda)F_1,\lambda E_2 +(1 -\lambda)F_2)
\end{equation}
holds.
\end{prop}

The case $F(t)=t^\alpha$ is the famous Lieb's concavity theorem \cite{Lieb}.

Our next aim is to compute 
\begin{equation}\label{E:deri1}
\frac{\pard^2}{\pard t \pard s}  S_F (D+tA,D+sB)\Big|_{t=s=0}.
\end{equation}
We shall use the formulas
$$
{d\over dt} h(D+tB)\Big|_{t=0} =Bh'(D)\quad (B \in \iM_D), \qquad
{d\over dt} h(D+t\im [D,X])\Big|_{t=0} =\im [h(D),X]\,,
$$
see \cite{qinfobook}.

\begin{lemma}
If $A,B \in \iM_D$, then the derivative (\ref{E:deri1}) equals 
$-F''(1)\Tr D^{-1}AB$.
\end{lemma}

\proof
It is enough to check the case $F(t)=t^n$. Then
\begin{equation}\label{E:deri2}
S_F (D+tA,D+sB)=\Tr (D+tA)^{1-n}(D+sB)^n
\end{equation}
and the derivative is $(1-n)n \Tr D^{-1}AB$. \qed

\begin{lemma}
If $A \in \iM_D$ and $B \in \iM_D^c$, then the derivative (\ref{E:deri1}) 
equals $0$.
\end{lemma}

\proof
We compute for $F(t)=t^n$ using (\ref{E:deri2}). If $B=[D,X]$, then we have
the derivative
$$
\Tr (1-n)D^{-n}A[D^n,X]=0
$$
and this gives the statement. \qed

\begin{lemma}
Let $X=X^* \in \iM$ and $F:\bbbr^+ \to \bbbr$ be a continously differentiable
function. Then
\begin{equation}\label{E:deri}
\frac{\pard^2}{\pard t \pard s}  S_F (D+t\im [D,X],D+s\im [D,X])\Big|_{t=s=0}
=2 F(1)\Tr DX^2-2 S^X_F (D,D).
\end{equation}
\end{lemma}

\proof
Since both sides are linear and continuous in $F$, we may assume that
$F(t)=t^n$. Derivation of formula (\ref{E:deri2}) gives 
$$
2 \Tr DX^2-2 \Tr XD^{1-n}XD^n
$$
and this is the stated result for the particular $F$. \qed

\section{Skew information}

The Wigner-Yanase-Dyson skew information is the quantity
$$
I_p(D,A):= -\frac{1}{2} \Tr [D^p, A][D^{1-p},A] \qquad (0 < p <1).
$$
Actually, the case $p=1/2$ is due to Wigner and Yanase \cite{WYD} 
and the extension was proposed by Dyson. The convexity of $I_p(D,A)$ in 
$A$ is a famous result of Lieb \cite{Lieb}

It was observed in \cite{PH} that the Wigner-Yanase-Dyson skew information
is connected to a monoton Riemannian metric (or Fisher information) which
corresponds to the function
$$
f_p(x)={p(1-p)}\,{ (x-1)^2 \over (x^p -1)
(x^{1-p}-1)}\, .
$$
It was proven in \cite{PH} that this is an operator monotone function, a 
generalization was obtained in \cite{H, Sa}.

Let $f$ be a standard function and $X=X^*\in \iM$. The quantity
$$
I_D^f(X):= \frac{f(0)}{2}\gamma_D^f(\im [D,X], \im [D,X] )
$$
was called {\bf skew information} in \cite{H} in this general setting.
Note that the parametrization in \cite{H} is by $c=1/f$ which is called
there Morozova-Chentsov function.
The skew information is nothing else but the Fisher information restricted
to $\iM_D^c$, but it is parametrized by the commutator. Skew information 
appears, for example, in uncertainty relations \cite{andai, GI, GII, 
GII:2007a, Kosaki, LZ1, LZ2}, see also Theorem \ref{T:uj}. In that 
application, the skew information is regarded as a bilinear form.

If $D=\Diag (\lambda_1, \dots, \lambda_n)$ is diagonal, then
$$
\gamma_D^f(\im [D,X], \im [D,X] )=\sum_{ij} 
\frac{(\lambda_i -\lambda_j)^2}{\lambda_j f(\lambda_i/\lambda_j)}|X_{ij}|^2.
$$
This implies that the identity
\begin{equation}\label{E:tilde2}
f(0)\gamma_D^f(\im [D,X], \im [D,X] )=
2\Cov_D(X,X)-2\QCov^{\tilde f}_D (X, X )
\end{equation}
holds if $\Tr DX=0$ and
\begin{equation}\label{E:tilde1}
\tilde{f}(x):=\frac{1}{2}\left((x+1)-(x-1)^2 \frac{f(0)}{f(x)}\right).
\end{equation}
Since the right-hand-sides of (\ref{E:deri}) and (\ref{E:tilde2}) are
the same  if $F=\tilde f$ we have

\begin{Thm}
Assume that $X=X^* \in \iM$ and $\Tr D X=0$. If $f$ is a standard function
such that $f(0)\ne 0$, then
$$
\frac{\pard^2}{\pard t \pard s}  S_F (D+t\im [D,X],D+s\im [D,X])\Big|_{t=s=0}
=f(0)\gamma_D^f(\im [D,X], \im [D,X] )
$$
for the standard function $F=\tilde f$.
\end{Thm}

The only remaing thing to show is that  if $f:\bbbr^+ \to \bbbr$ is a 
standard function, then $\tilde{f}$ is standard as well. This result 
appeared in \cite{3} and the proof there is not easy, even matrix convexity
of functions of two variables is used. Here we give a rather elementary 
proof based on the fact that $1/f  \mapsto \tilde{f}$ is linear and 
on the canonical decomposition in Theorem \ref{T:hansen}. 

\begin{lemma}
Let $\lambda\ge0$ and $f:\bbbr^+ \to \bbbr$ be a function such that
$$
\frac{1}{f(x)}:=
\frac{1+\lambda}{2}\left(\frac{1}{x+\lambda}+\frac{1}{1+x\lambda}\right)
=g_\lambda(x).
$$
Then the function $\widetilde{f}:\bbbr^+ \to \bbbr$ defined in (\ref{E:tilde1})
is an operator monotone standard function.
\end{lemma}

\proof
From the definitions we obtain
\[
\widetilde{f}(x)=\frac{x(x\lambda^{2}+\lambda^{2}+2\lambda+2x\lambda+x+1)}
{2(x+\lambda)(1+x\lambda)}
\]
and
\[
\widetilde{f}'(x)=\frac{\lambda+2x\lambda+2\lambda^{2}+x^{2}\lambda+
4x\lambda^{2}+2\lambda^{3}x+x^{2}+\lambda^{3}x^{2}+\lambda^{3}+
\lambda^{4}x^{2}}{2(x+\lambda)^{2}(1+x\lambda)^{2}}.
\]
Hence $\widetilde{f}(0)=0$ and $\widetilde{f}:\bbbr^+ \to \bbbr^+$. So it is 
enough to prove that the holomorphic extension of $\widetilde{f}$ 
to the complex upper half-plane maps the upper half-plane into itself,
see \cite{Bh}.

Let $a,b\in\bbbr$, $b>0$. Then we have
\begin{eqnarray*}
\Ima \widetilde{f}(a+ib)&=&
\frac{b}{2((a+\lambda)^{2}+b^{2})((1+\lambda a)^{2}+\lambda^{2}b^{2})}
\cr && \times
(\lambda+2\lambda^{2}+\lambda^{3}+b^{2}+a^{2}+2\lambda a+\lambda a^{2}
\cr && \quad
+4\lambda^{2}a+\lambda b^{2}+\lambda^{3}a^{2}+\lambda^{3}b^{2}+2\lambda^{3}a
+\lambda^{4}a^{2}+\lambda^{4}b^{2}).
\end{eqnarray*}
Here
\begin{eqnarray*}
\lambda+2\lambda^{2}+\lambda^{3}+b^{2}+a^{2}+2\lambda a+\lambda a^{2}+
4\lambda^{2}a+\lambda b^{2}+\lambda^{3}a^{2}+\lambda^{3}b^{2}+2\lambda^{3}a+
\lambda^{4}a^{2}+\lambda^{4}b^{2}\cr
=(1+\lambda+\lambda^{4}+\lambda^{3})a^{2}+(4\lambda^{2}+2\lambda+2\lambda^{3})a
+\lambda+2\lambda^{2}+\lambda^{3}+b^{2}(1+\lambda+\lambda^{4}+\lambda^{3}).
\end{eqnarray*}
The function $g:\bbbr \to \bbbr$ 
\[
g(a):=(1+\lambda+\lambda^{4}+\lambda^{3})a^{2}+(4\lambda^{2}+2\lambda+2\lambda^{3})a+\lambda+2\lambda^{2}+\lambda^{3}
\]
has a minimum value at 
\[
a(\lambda)=-\frac{4(\lambda^{2}+2\lambda+2\lambda^{3})}{2(1+\lambda+\lambda^{4}+\lambda^{3})}
\]
 and 
\[
g(a(\lambda))=\frac{(\lambda^{2}-1)^{2}\lambda}{(\lambda-1/2)^{2}+3/4}\geq0.
\]
Therefore the upper half-plane is mapped into itself. The properties
$x\widetilde{f}(x^{-1})=\widetilde{f}(x)$ and $\widetilde{f}(1)=1$ 
are obvious. \qed

The uncertainty relation recently obtained is the following \cite{GHP}.

\begin{prop}\label{T:uj}
Assume that $f,g:\bbbr^+\to \bbbr$ are standard functions 
and $D$ is a positive definite matrix. Then for self-adjoint matrices 
$A_1,A_2, \dots ,A_m$ the determinant 
inequality
\begin{eqnarray*}
&&\Det\biggl( \left[\QCov_{D}^g(A_i,A_j)\right]_{i,j=1}^m \biggr) 
\cr && \qquad\ge
\Det\biggl( \left[ f(0)g(0)
\big(\Cov_D(A_i,A_j)-\QCov^{\tilde f}_D (A_i, A_j )\big)
\right]_{i,j=1}^m  \biggr)
\end{eqnarray*}
holds.  
\end{prop}

Note that the right-hand-side contains skew informations, 
cf. (\ref{E:tilde2}).

\section{The setting of von Neumann algebras}

Let $\iM$ be a von Neumann algebra. Assume that it is in standard form,
it acts on a Hilbert space $\iH$, $\iP\subset \iH$ is the positive cone 
and $J:\iH \to \iH$ is the modular conjugation \cite{Uffe, OP, St}. 
Let $\ffi$ and $\omega$ 
be normal states with representing vectors $\Phi$ and $\Omega$ in the 
positive cone. For the sake of simplicity, assume that $\ffi$ and 
$\omega$ are faithful. This means that $\Phi$ and $\Omega$ are cyclic 
and separating vectors. The closure of the unbounded operator 
$A\Omega \mapsto A^*\Phi$ has a polar decomposition $J\Delta(\ffi, 
\omega)^{1/2}$ and $\Delta(\ffi, \omega)$ is called relative modular 
operator. $A\Omega $ is in the domain of $\Delta(\ffi, \omega)^{1/2}$ 
for every $A \in \iM$.

For $A \in \iM$ and $F:\bbbr^+ \to \bbbr$, the quasi-entropy
\begin{equation}\label{E:quasi2}
S^A_F (\omega,\ffi):= \< A \Omega, F(\Delta(\ffi, \omega))A\Omega\> 
\end{equation}
was introduced in \cite{PD26}, see also Chapter 7 in \cite{OP}. 
(The right-hand-side can be understood via the spectral decomposition
of the positive operator $\Delta(\ffi, \omega)$.) For $F(t)=-\log t$ and 
$A=I$ the relative entropy of Araki is obtained \cite{Araki} and
this was the motivation of the generalization.

\begin{Thm}\label{T:quasimon}
Assume that $F:\bbbr^+ \to \bbbr$ is an operator monotone function
with $F(0)\ge 0$ and $\alpha:\iM_0 \to \iM$ is a unital normal Schwarz 
mapping. Then
\begin{equation}\label{E:quasimon2}
S^A_F (\omega\circ \alpha,\ffi\circ \alpha) \ge
S^{\alpha(A)}_F (\omega,\ffi)
\end{equation}
holds for $A \in \iM_0$ and for normal states $\omega$ and $\ffi$
of the von Neumann algebra $\iM$.
\end{Thm}

We sketch the proof based on inequalities for operator monotone and
operator concave functions. (The details are clarified in \cite{PD26}.) 
First note that
$$
S^A_{F+c} (\omega\circ \alpha,\ffi\circ \alpha)=
S^A_F (\omega\circ \alpha,\ffi\circ \alpha)+c\,\omega(\alpha(A^*A))
$$
and 
$$
S^{\alpha(A)}_{F+c} (\omega,\ffi)=S^{\alpha(A)}_F (\omega,\ffi)+
c\,\omega(\alpha(A)^*\alpha(A))
$$
for a positive constant $c$. Due to the Schwarz inequalty (\ref{E:Sch}),
we may assume that $F(0)=0$.

Let $\Omega_0$ be the representing vector for $\omega \circ \alpha$
and $\Delta:=\Delta(\ffi, \omega)$, $\Delta_0:=\Delta(\ffi\circ \alpha, 
\omega \circ \alpha)$. The operator
\begin{equation}
Vx\Omega_0=\alpha(x)\Omega \qquad (x \in \iM_0)
\end{equation}
is a contraction:
$$
\Vert \alpha(x)\Omega \Vert^2 = \omega (\alpha(x)^* \alpha(x))\le
\omega (\alpha(x^*x)= \Vert x \Omega_0 \Vert^2
$$
since the Schwarz inequality is applicable to $\alpha$. A similar simple
computation gives that
\begin{equation}
V^*\Delta V \le \Delta_0\,.
\end{equation}

Since $F$ is operator monotone, we have $F(\Delta_0)\ge F(V^*\Delta V)$. 
Recall that $F$ is operator concave, therefore $F(V^*\Delta V) \ge 
V^*F(\Delta)V$ and we conclude
\begin{equation}
F(\Delta_0) \ge V^*F(\Delta)V\,.
\end{equation}
Application to the vector $A \Omega_0$ gives the inequality.

The natural extension of the covariance (from probability theory) is
\begin{equation}
\QCov^f_\omega(A,B)=\< \sqrt{f(\Delta(\omega, \omega))}A \Omega, 
\sqrt{f(\Delta(\omega, \omega))}B\Omega\>
-\overline{\omega(A)}\omega(B),
\end{equation}
where $\Delta(\omega, \omega)$ is actually the modular operator.
Motivated by the application, we always assume that the function $f$ is 
standard. For such a function $f$, the inequalities
$$
\frac{2x}{x+1} \le f(x) \le \frac{1+x}{2}
$$ 
holds. Therefore $A\Omega$ is in the domain of $\sqrt{f(\Delta(\omega, 
\omega))}$ and the covariance $\QCov^f_\omega(A,B)$ is a well-defined 
sesquilinear form.

For a standard function $f:\bbbr^+ \to \bbbr^+$ and for a normal
unital Schwarz mapping $\beta: \iN \to \iM$ the inequality
\begin{equation}
\QCov^f_\omega(\beta(X),\beta(X))
\le
\QCov^f_{\omega\circ \beta}(X, X) \qquad (X \in \iN)
\end{equation}
is a particular case of Theorem \ref{T:quasimon} and it is the monotonicity 
of the generalized covariance under coarse-graining \cite{PD22}.

Following \cite{H}, the skew information (as a bilinear form) can be defined as
\begin{equation}\label{E:skew_vN}
I_\omega^f(X, Y):= 
\Cov_\omega(X,Y)-\QCov^{\tilde f}_\omega (X, Y)
\end{equation}
if $\omega(X)=\omega(Y)=0$. (Then $I_\omega^f(X)=I_\omega^f(X,X)$.)

\begin{Thm}\label{T:vegso}
Assume that $f,g:\bbbr^+\to \bbbr$ are standard functions 
and $\omega$ is a faithful normal state on a von Neumann algebra $\iM$.
Let $A_1,A_2, \dots ,A_m\in \iM$ be self-adjoint operators such that
$\omega(A_1)=\omega(A_2)= \dots =\omega(A_m)=0$. Then the determinant 
inequality
\begin{eqnarray}
&&\Det\biggl( \left[\QCov_{D}^g(A_i,A_j)\right]_{i,j=1}^m \biggr) 
\ge
\Det\biggl( \left[ 2g(0)I_\omega^f(A_i, A_j)\right]_{i,j=1}^m  \biggr)
\end{eqnarray}
holds.  
\end{Thm}

\proof
Let $E(\pont)$ be the spectral measure of $\Delta(\omega, \omega)$. Then
for $m=1$ the inequality is
$$
\int g(\lambda) \,d\mu(\lambda)
\le g(0)\left(
\int \frac{1+\lambda}{2} \,d \mu(\lambda)
-\int {\tilde f}(\lambda) \,d \mu(\lambda)\right),
$$
where $d\mu(\lambda)=d \< A\Omega, E(\lambda)A\Omega\>$. Since the inequality
\begin{equation}\label{E:gibi}
f(x)g(x)\ge f(0)g(0) (x-1)^2
\end{equation}
holds for standard functions \cite{GHP}, we
have
$$
g(\lambda)\ge g(0)\left(\frac{1+\lambda}{2}-f(0){\tilde f}(\lambda)\right)
$$
and this implies the integral inequality.
 
Consider the finite dimensional subspace $\iN$ generated by the operators
$A_1,A_2, \dots ,A_m$. On $\iN$ we have the inner products
$$
\bal A, B \jobb :=\Cov_{\omega}^g(A,B)
$$
and
$$
\< A,B\>:= 2 g(0)I_\omega^f(A, B).
$$
Since $\< A,A\> \le \bal A, A \jobb$, the determinant inequality 
holds (see Lemma 2 in \cite{GHP}).\qed

This theorem is interpreted as quantum uncertainty principle
\cite{andai, GII, 3, Kosaki}. In the earlier works the function $g$ from the
left-hand-side was $(x+1)/2$ and the proofs were more complicated. 
The general $g$ appeared in \cite{GHP}.


\begin{thebibliography}{99}

\bibitem{andai}
A. Andai, Uncertainty principle with quantum Fisher information,
to be published in J. Math. Phys. 

\bibitem{Araki}
H. Araki, Relative entropy of states of von Neumann algebras,  
Publ. Res. Inst. Math. Sci.  {\bf 11}(1975/76), 809--833.

\bibitem{Bh}
R. Bhatia, {\it Matrix Analysis}, Springer, 1997.

\bibitem{Csi}
I. Csisz{\'a}r, Information type measure of difference of
probability distributions and indirect observations, Studia 
Sci. Math. Hungar. {\bf 2}(1967), 299--318.

\bibitem{GI}  
P. Gibilisco and T. Isola,  Uncertainty principle and
quantum Fisher information, Ann.  Inst.  Stat.  Math, {\bf 59}
(2007), 147--159.

\bibitem{GII}
P. Gibilisco, D. Imparato and T. Isola,  A volume inequality for quantum 
Fisher information and the uncertainty principle, to be published in 
J. Statist.  Phys.

\bibitem{GII:2007a}
P. Gibilisco, D. Imparato and T. Isola, A Robertson-type uncertainty 
principle and quantum Fisher information, to be published in  Lin. Alg. Appl.

\bibitem{3}
P. Gibilisco, D. Imparato and T. Isola, Uncertainty principle and 
quantum Fisher information II, J. Math. Phys. {\bf 48}(2007), 072109, 
arXiv:math-ph/0701062v3.

\bibitem{GHP}
P. Gibilisco, F. Hiai and D. Petz, Quantum covariance, quantum Fisher 
information and the uncertainty principle, preprint, 2007.

\bibitem{Uffe}
U. Haagerup, The standard form of von Neumann algebras,  Math. Scand.  
{\bf 37}(1975), 271--283.

\bibitem{H1}
F. Hansen, Characterizations of symmetric monotone metrics on the the state
space of quantum systems, Quantum Inf. Comput., {\bf 6}(2006), 597--605.

\bibitem{H}
F. Hansen, Metric adjusted skew information, arXiv:math-ph/0607049v3, 2006.

\bibitem{Kosaki} 
H. Kosaki, Matrix trace inequality related to uncertainty principle,
Internat. J. Math. {\bf 16}(2005), 629--645.

\bibitem{Ando}
F. Kubo and T. Ando, Means of positive linear operators, 
Math. Ann. {\bf 246}(1980), 205--224.

\bibitem{Lieb}
E. H. Lieb, Convex trace functions and the
Wigner-Yanase-Dyson conjecture, Advances in Math. {\bf 11}(1973), 267--288.

\bibitem{LV} 
F. Liese and I. Vajda, On divergences and informations in statistics and 
information theory, IEEE Trans. Inform. Theory {\bf 52}(2006), 4394-4412.

\bibitem{LZ1} 
S. Luo and Z. Zhang,
An informational characterization of Schr\"odinger's uncertainty relations,
J. Stat. Phys. {\bf 114}, 1557--1576 (2004).

\bibitem{LZ2} 
S. Luo and Q. Zhang, On skew information,
IEEE Trans. Inform. Theory, {\bf 50}(2004), 1778--1782.

\bibitem{LZ3} S. Luo and Q. Zhang.
Correction to On skew information.
IEEE Trans. Inform. Theory {\bf 51}, 4432 (2005).

\bibitem{OP}
M. Ohya and D. Petz,
{\it Quantum Entropy and Its Use},
Springer-Verlag, Heidelberg, 1993. Second  edition 2004.

\bibitem{PD26}
D. Petz,
Quasi-entropies for states of a von Neumann algebra,
Publ. RIMS. Kyoto Univ. {\bf 21}(1985), 781--800.

\bibitem{PD32} 
D. Petz, Quasi-entropies for finite quantum systems,
Rep. Math. Phys., {\bf 23}(1986), 57-65.

\bibitem{PD2}
D.~Petz,   Monotone metrics on matrix spaces. Linear
Algebra Appl. {\bf 244}(1996), 81--96.

\bibitem{PH}
D. Petz and H. Hasegawa, On the Riemannian metric of $\alpha$-entropies 
of density matrices, Lett. Math. Phys. {\bf 38}(1996), 221--225

\bibitem{PD22}
D.~Petz,  Covariance and Fisher information in quantum mechanics.
J. Phys. A: Math. Gen. {\bf 35}(2003), 79--91.

\bibitem{qinfobook}
D. Petz, {\it Quantum Information Theory and Quantum Statistics},
Springer-Verlag, Berlin and Heidelberg, 2007.

\bibitem{St}
\c S. Str\u atil\u a, {\it Modular theory in operator algebras},
Abacus Press, Tunbridge Wells, 1981. 

\bibitem{Sa}
V. E. S. Szab\'o, A class of matrix monotone functions.  
Linear Algebra Appl.  {\bf 420}(2007), 79--85. 

\bibitem{WYD}
E. P. Wigner and M. M. Yanase,  Information content of distributions.
Proc. Nat. Acad. Sci. USA {\bf 49}(1963), 910--918.

\end{thebibliography}
\end{document}